\documentstyle{article}
\newcommand{\be}{\begin{equation}}
\newcommand{\ee}{\end{equation}}
\newcommand{\ba}{\begin{eqnarray}}
\newcommand{\ea}{\end{eqnarray}}
\newcommand{\baa}{\begin{eqnarray*}}
\newcommand{\eaa}{\end{eqnarray*}}
\newcommand{\bb}{\begin{thebibliography}}
\newcommand{\eb}{\end{thebibliography}}

\newcommand{\bi}[1]{\bibitem{#1}}
\newcommand{\lab}[1]{\label{#1}}
\newcommand{\re}[1]{(\ref{#1})}

\newcounter{my}
\newcommand{\he}%
   {\stepcounter{equation}\setcounter{my}%
   {\value{equation}}\setcounter{equation}0%
   }%
\newcommand{\she}%
   {\setcounter{equation}{\value{my}}%
    }%

\newtheorem{lem}{Lemma}
\newtheorem{tm}{Theorem}

\begin{document}
\begin{center}
{\Large {\bf Generalized Bochner theorem: characterization of the
Askey-Wilson polynomials}}  \vspace{5mm}

\medskip

{\bf Luc Vinet}

\medskip

{\it Universit\'e de Montr\'eal PO Box 6128, Station Centre-ville
Montréal QC H3C 3J7. e-mail: luc.vinet@umontreal.ca}
\medskip

and

\medskip

{\bf  Alexei Zhedanov}

\medskip

{\it Donetsk Institute for Physics and Technology, Donetsk 83114,
Ukraine. e-mail: zhedanov@kinetic.ac.donetsk.ua}

\bigskip

\medskip

\end{center}

\begin{abstract}
Assume that there is a set of monic polynomials $P_n(z)$
satisfying the second-order difference equation
$$
A(s) P_n(z(s+1)) +  B(s) P_n(z(s)) + C(s) P_n(z(s-1)) = \lambda_n
P_n(z(s)), \; n=0,1,2,\dots, N$$ where $z(s), A(s), B(s), C(s)$
are some functions of the discrete argument $s$ and $N$ may be
either finite or infinite. The
irreducibility condition $A(s-1)C(s) \ne 0$ is assumed for all
admissible values of $s$. In the finite case we assume that there
are $N+1$ distinct grid points $z(s), \: s=0,1,\dots, N$ such that
$z(i) \ne z(j), \: i \ne j$. If $N=\infty$ we assume that the grid
$z(s)$ has infinitely many different values for different values
of $s$. In both finite and infinite cases we assume also that the
problem is non-degenerate, i.e. $\lambda_n \ne \lambda_m, \; n \ne
m$. Then we show that necessarily: (i) the grid $z(s)$ is at most
quadratic or q-quadratic in $s$; (ii) corresponding polynomials
$P_n(z)$ are at most the Askey-Wilson polynomials corresponding to
the grid $z(s)$.  This result can be considered as generalizing of
the Bochner theorem (characterizing the ordinary classical
polynomials) to generic case of arbitrary difference operator on
arbitrary grids.

\vspace{15mm}



\end{abstract}

\newpage

\section{Introduction}
General orthogonal polynomials (OP) $P_n(x)$ can be characterized
by the 3-term recurrence relation \cite{Chi} \be P_{n+1}(x) + b_n
P_n(x) + u_nP_n(x) = xP_n(x) \lab{3term} \ee with initial
conditions $P_0=1, P_1=x-b_0$

The polynomials $P_n(x)$ are {\it monic} polynomials, i.e. $P_n(x)
= x^n + O(x^{n-1})$.

It is well known \cite{Al} that all polynomials solutions $P_n(x)$
of the second-order differential equation \be \sigma(x) P_n''(x) +
\tau(x) P_n'(x) = \lambda_n P_n(x) \lab{hyp} \ee are classical
orthogonal polynomials (COP), i.e. Jacobi, Laguerre, Hermite and
Bessel polynomials.  In \re{hyp} it appears that $\sigma(x)$ and
$\tau(x)$ are polynomials such that $deg(\sigma(x) \le 2, \;
deg(\tau(x) = 1$.  This result is known as the Bochner theorem
\cite{Bo}.

It is natural to consider generalization of the Bochner theorem
replacing the second-order differential operator with the
second-order difference operator. Namely we are seeking polynomial
solutions $P_n(z)$ of the problem \be A(s) P_n(z(s+1)) +  B(s)
P_n(z(s)) + C(s) P_n(z(s-1)) = \lambda_n P_n(z(s)), \;
n=0,1,2,\dots, N \lab{gen_B} \ee where $z(s), A(s), B(s), C(s)$
are some functions of the discrete argument $s$ and $N$ may be
either finite or infinite. The irreducibility condition
$A(s-1)C(s) \ne 0$ is assumed for all admissible values of $s$. In
the finite case we assume that there are $N+1$ distinct grid point
$z(s), \: s=0,1,\dots, N$ such that $z(i) \ne z(j), \: i \ne j$.
If $N=\infty$ we assume that the grid $z(s)$ has infinitely many
different values for different values of $s$. In both finite and
infinite cases we assume also that the problem is non-degenerate,
i.e. $\lambda_n \ne \lambda_m, \; n \ne m$. We assume also that
there are polynomial solutions of all degrees $n=0,1,\dots, N$
(i.e. we assume that the polynomial $P_n(x)$ always has exact
degree $n$ for all $n=0,1,\dots, N$.)

Askey and Wilson  \cite{AW} discovered orthogonal polynomials (the
Askey-Wilson polynomials, or briefly, AWP) which satisfy equation
\re{gen_B} for quadratic $z(s) = as^2 + bs +c$ or q-quadratic grid
$z(s) = aq^s + b q^{-s} + c$, where $q$ is some parameter such
that $|q| \ne 1$. Finite-dimensional case (i.e. when there exists
only $N$ mutually orthogonal polynomials $n=0,1,\dots,N-1$)
corresponds to the so-called q-Racah polynomials \cite{KS}.

In \cite{GH}  it was shown that the only OP satisfying \re{gen_B}
for AW-grids are the AWP. Leonard showed \cite{Leo} that in the
finite-dimensional case the only OP satisfying  \re{gen_B} are the q-Racah polynomials.
 For further development of the Leonard result and its
new algebraic interpretation see, e.g. \cite{Terw}. In
\cite{Ismail} Ismail obtained more strong result : he showed that
all polynomial (i.e. not necessarily orthogonal, {\it ab initio})
solutions of the equation \re{gen_B} for the AW-grid are AWP. In the finite-dimensional case 
Terwilliger obtained the result that the AW-grid is the most general for polynomials 
satisfying \re{gen_B}.  

So far, the open problem was: in the infinite-dimensional case 
characterize all possible grids
$z(s)$ for which polynomial solutions of the equaiton \re{gen_B}
are obtained. In this paper we solve this problem and show that
there are no grids more general than AW-grids. Hence  all
polynomial solutions for \re{gen_B} {\it should be} orthogonal
Askey-Wilson polynomials. Although for the finite-dimensional case the problem was effectively solved by Terwilliger in \cite{Terw1},
we present here the finite-dimensional version of the generalized Bochner theorem as well. 
The main reason is that our method
of proof is essentially different and deals directly with difference equation \re{gen_B} for polynomials whereas in the Terwilliger paper \cite{Terw1} another (a purely algebraic) approach is presented.

\section{Finite-dimensional case}
\setcounter{equation}{0} In this section we show that if $N$ is
finite then the problem is essentially equivalent to the Leonard
theorem \cite{Leo}.

Indeed, consider $(N+1) \times (N+1)$ tri-diagonal matrix $J$
which acts on a basis $e_k, \: k=0,1,\dots,N$ by \be J e_k =
C(k+1) e_{k+1} + B(k) e_k + A(k-1) e_{k-1} \lab{def_J} \ee It is
assumed that $C(N+1)=A(-1)=0$ which means merely that that the
matrix $J$ acts in linear space of dimension $N+1$. We will assume
the nondegeneracy condition: \be C(i)A(i-1) \ne 0 , \; i=1,2\dots
N \lab{ndeg_AC} \ee

Find the eigenvectors $v^{(k)}, k=0,1,\dots, N$ of the matrix $J$,
i.e. $$J v^{(k)} = \lambda_k v^{(k)}$$ with some eigenvalues
$\lambda_k$. We assume that all eigenvalues are distinct:
$\lambda_i \ne \lambda_j$ if $i \ne j$. Then all vectors $v^{(k)},
\: k=0,1,\dots, N$ are independent and we have \be
v^{(k)}=\sum_{s=0}^N v_{ks} e_s, \lab{expans_v} \ee where $v_{ks},
s=0,1,\dots,N$ are components of the vector $v^{(k)}$ in the basis
$e_s$. For them we have relation \be A(s) v_{k,s+1} + B(s) v_{ks}
+ C(s) v_{k,s-1} = \lambda_k v_{ks}. \lab{EVv} \ee  Now we can
identify components $v_{ks}$ with $P_k(z_s)$, i.e. we merely put
$v_{ks} = P_k(z_s)$ for all values $k,s=0,1,\dots, N$. Then
difference equation \re{gen_B} coincides with \re{EVv}.

Consider transposed Jacobi matrix $J^*$ defined as \be J^* e_k =
A(k)e_{k+1} + B(k) e_k + C(k) e_{k-1} \lab{conj_J} \ee and
corresponding eigenvalue vectors $v^{*(k)}$: \be J^* v^{*(k)} =
\lambda_k \: v^{*(k)}, \quad k=0,1,2,\dots, N \lab{conj_Jv} \ee
Vectors $v^{*(k)}$ can be expanded in terms of the same basis
$e_s$: \be v^{*(k)} = \sum_{s=0}^N v^*_{ks} e_s \lab{v_conj} \ee

From elementary linear algebra it is known that in nondegenerated
case (i.e. if $\lambda_i \ne \lambda_j$ for $i \ne j$) the vectors
$v^k$ and $v^{*(j)}$ are biorthogonal: \be (v^k,v^{*(j)}) \equiv
\sum_{s=0}^N v_{ks} v^*_{js} = 0 \quad \mbox{if} \quad k \ne j
\lab{bi_vv} \ee Introduce now the diagonal matrix $M$ wich acts on
basis $e_s$ as \be M e_s = \mu_s e_s, \quad s=0,1,2,\dots, N
\lab{M_def} \ee where \be \mu_s =\frac{A(0)A(1) \dots
A(s-1)}{C(1)C(2) \dots C(s)}, \; s=1,2,\dots,N, \quad \mu_0 =1
\lab{mu_def} \ee Note that all $\mu_s$ are well defined due to
nondegeneracy condition \re{ndeg_AC}.

It is elementary verified that \be J^* = M^{-1} J M,  \lab{JMJ}
\ee and hence \be v^{*(k)} = M^{-1} v^{(k)}, \quad k=0,1,\dots, N
\lab{vMv} \ee (inverse matrix $M^{-1}$ exists due to nondegeneracy
condition \re{ndeg_AC}). Relation \re{vMv} allows one to rewrite
biorthogonality condition \re{bi_vv} in the form \be \sum_{s=0}^N
w_s v_{ks} v_{js} =0, \quad \mbox{if} \quad k \ne j, \lab{ort_vv}
\ee where \be w_s = 1/\mu_s =  \prod_{i=1}^s \frac{C(i)}{A(i-1)}
\lab{ws} \ee In terms of polynomials $P_n(x)$ this relation
becomes \be \sum_{s=0}^N w_s P_k(z(s)) P_j(z(s)) =0, \quad
\mbox{if} \quad k \ne j \lab{ort_P} \ee But relation \re{ort_P}
means that $P_n(x)$ are polynomials which are orthogonal on a
finite distinct set of points $z(s), s=0,1,\dots,N$ with discrete
weights $w_s \ne 0$. By general elementary theorems concerning
orthogonal polynomials \cite{Chi} this means that polynomials
$P_n(x)$ should satisfy a three-term recurrence relation \be
P_{n+1}(x) + b_n P_n(x) + u_n P_{n-1}(x) = xP_n(x), \quad
n=0,1,\dots, N \lab{3_term_f} \ee The roots $x_s, \: s=0,1,\dots,
N$ of the polynomial $P_{N+1}(x)$ coincide with spectral points:
$$
z(s) = x_s, \quad s=0,1,\dots, N
$$
Thus we proved that (under some nondegeneration conditions)
polynomials $P_n(x)$ satisfying relation \re{gen_B} on a grid
$z(s)$ for finite $N$ are orthogonal with respect to discrete
weight function \re{ort_P} and satisfy three-term recurrence
relation \re{3_term_f}.

Now we are ready to relate our results with Leonard's approach to
dual orthogonal polynomials \cite{Leo}.

Recall relation between nondegenerated Jacobi matrices and
orthogonal polynomials (see, e.g. \cite{Chi}). Let $K$ be an
arbitrary Jacobi matrix of dimension $N+1 \times N+1$. In some
finite-dimensional basis $d_n$ it can be presented as \be K d_n =
\alpha_n d_{n+1} + \beta_n d_n + \gamma_n d_{n-1} \lab{Jac_K} \ee
with some (complex) coefficients with nondegeneracy property \be
\prod_{i=1}^N \gamma_i \alpha_{i-1} \ne 0 \lab{n_deg_K} \ee
Construct eigenvectors $\pi^{(k)}$ of the matrix $K$: \be K
\pi^{(k)} = z_k \pi^{(k)}, \quad k=0,1,\dots, N \lab{pi_K} \ee We
assume that all spectral points $z_k, \: k=0,1,\dots,N$ are
distinct: $z_k \ne z_j$ for $k \ne j$. Expand eigenvectors
$\pi^{(k)}$ in terms of basis $d_n$:
$$
\pi^{(k)}= \sum_{s=0}^N \pi_{ks} d_s
$$
with some coefficients $\pi_{ks}$. For these coefficients we have
from \re{pi_K} the recurrence relation \be \gamma_{s+1}
\pi_{k,s+1} + \beta_s \pi_{ks} + \alpha_{s-1} \pi_{k,s-1} = z_k
\pi_{ks}, \quad k,s=0,1, \dots, N \lab{rec_pi} \ee It is assumed
that $\alpha_{-1}=\gamma_{N+1}=0$. Then, for each value $k$,
starting from $\pi_{k0}$ we can find recursively all further
coefficients $\pi_{k1},\pi_{k,2},\dots, \pi_{k,N}$. We can always
normalize $\pi_{k0}=1, \: k=0,1,\dots,N$. Then it is clear from
\re{rec_pi} that $\pi_{ks}$ is a polynomial of degree $s$ in
argument $z_k$.

Introduce polynomials $T_n(x)$ satisfying three-term recurrence
relation \be \gamma_{n+1} T_{n+1}(x) + \beta_n T_n(x) +
\alpha_{n-1} T_{n-1}(x) = x T_n(x) \lab{rec_T} \ee with initial
conditions $\alpha_{-1}P_{-1}=0, \: P_0(x)=1$. Then relation
\re{rec_T} defines $n$-degree polynomials $T_n(x)=\kappa_n x^n +
O(x^{n-1})$ with the leading coefficient
$$
\kappa_n = \frac{1}{\gamma_1 \gamma_2 \dots \gamma_n}
$$
(this leading coefficient is well defined and nonzero do to
nondegeneracy condition \re{n_deg_K}). From general theory of
orthogonal polynomials it follows that polynomials $T_n(x)$ are
orthogonal on a finite set of points $x_k$ \cite{Chi} \be
\sum_{k=0}^N \sigma_k T_n(x_k)T_m(x_k) = 0, \quad n \ne m
\lab{ort_T} \ee where $x_k$ are roots of the polynomial
$T_{N+1}(x)$.

We can thus associate orthogonal polynomials $T_n(x)$ with
expansion coefficients of eigenvectors of the arbitrary
nondegenerated Jacobi matrix $K$: \be T_s(z_k) = \pi_{ks}
\lab{T_pi} \ee

Return to our polynomials $P_n(x)$ satisfying relation \re{gen_B}.
We showed that these polynomials are orthogonal and correspond to
the Jacobi matrix $K$ whose matrix coefficients can be restored
from recurrence relation \re{3_term_f}: $\gamma_n =1, \beta_n =
b_n, \alpha_n = u_{n+1}$. On the other hand, we have the Jacobi
matrix $J$ defined by \re{def_J}. By just described recipe, we can
associate with this Jacobi matrix corresponding orthogonal
polynomials $Y_n(x)$. These polynomials satisfy 3-term recurrence
relation \be A(n) Y_{n+1}(x) + B(n) Y_n(x) + C(n) Y_{n-1}(x) =
xY_n(x) \lab{rec_Y1} \ee Now it is seen that polynomials $P_n(x)$
and $Y_n(x)$ are related as \be P_n(z(s)) = Y_s(\lambda_n)
\lab{duality} \ee We thus have a duality property coinciding with
that introduced by Leonard \cite{Leo}: there are two systems of
finite orthogonal polynomials and two finite sequences $z(s)$ and
$\lambda_n$ such relation \re{duality} holds. Our nondegeneracy
conditions: all $z(s)$ and $\lambda_n$ are distinct and matrices
$J,T$ are nondegenerated coincide with similar conditions in the
Leonard paper. Hence we can conclude
\begin{tm}
Under nondegeneracy conditions the finite-dimensional case of
relation \re{gen_B} generates at most finite Askey-Wilson
orthogonal polynomials (Racah and q-Racah polynomials in other
terms).
\end{tm}

\section{Infinite-dimensional case. Reducing to a more simple
problem} \setcounter{equation}{0} In this section we start to
analyze the infinite-dimensional case. We first derive some
restrictions upon the coefficients $A(s),B(s),C(s)$.

In what follows we will assume that polynomial solutions $P_n(z)$
of the equation \re{gen_B} are monic, i.e. $P_n(z)=z^n +
O(z^{n-1})$. This is not restriction of our problem, because it is
possible to divide all terms in equation \re{gen_B} by a (nonzero)
leading coefficient of the polynomial $P_n(z)$.

First of all we observe that eigenvalues $\lambda_n$ can be
shifted by an arbitrary constant $\lambda_n \to \lambda_n +
const$. Such shift leads to adding a constant to the coefficient
$B(s)$. Using this observation we always can choose $\lambda_n$ in
such a way that \be \lambda_0=0 \lab{l0} \ee In what follows we
will assume that condition \re{l0} is fulfilled. We will also
assume that the eigenvalue problem \re{gen_B} is non-degenerate,
i.e. \be \lambda_n \ne \lambda_m, \quad n \ne m \lab{ndegn} \ee
The grid $z(s)$ is also assumed to be non-degenerate, i.e. \be
z(s_1) \ne z(s_2), \quad s_1 \ne s_2 \lab{ndegs} \ee Parameter $s$
takes infinite number of integer values: $s=s_0,s_0+1,s_0+2,
\dots$ where $s_0$ is either finite or $s_0=-\infty$. In the first
case we deal with semi-infinite grid $z_s$, whereas in the second
case we have the grid which is infinite in both directions.

Taking the case $n=0$ in \re{gen_B} we see that $A(s) +B(s)+
C(s)=0$. Hence we can rewrite equation \re{gen_B} in the form \be
A(s) \Delta P_n(z(s)) - C(s) \nabla P_n(z(s)) = \lambda_n
P_n(z(s)), \lab{2form} \ee where we use the standard notation
\cite{NSU}
$$
\Delta F(s) = F(s+1) - F(s), \quad  \nabla F(s) = F(s) - F(s-1)
$$
for any function $F(s)$ of the argument $s$.

Assume that polynomials $P_n(z)$ have the expansion
$$
P_n(z) = z^n + \sum_{i=0} \xi_{ni} z^i
$$
with some coefficients $\xi_{ni}$. Then for $n=1$ we get from
\re{2form}
$$
A(s) \Delta z(s) - C(s) \nabla z(s) = \lambda_1 Q_1(z(s)),
$$
where $Q_1(z) = z + \xi_{10}= P_1(z)$.  By induction, it can be
easily shown that \be A(s) \Delta z^n(s) - C(s) \nabla z^n(s) =
\lambda_n Q_n(z(s)), \quad n=0,1,2,\dots, \lab{eq_Q} \ee where
$Q_n(z)$ is a monic polynomial of degree $n$.

Vice versa, assume that property \re{eq_Q} holds for some
$z(s),A(s),C(s)$ with $Q_n(x)$ being a set of monic polynomials in
$x$ of degree $n$. Then there exists a set of monic polynomials
$P_n(x)$ satisfying equation \re{2form}. This statement is almost
obvious and follows from the observation that on the given grid
$z(s)$ and for {\it any} monic $n$-th degree polynomial $T_n(x)$
the expression $A(s) \Delta T_n(z(s)) - C(s) \nabla T_n(z(s))$ is
again a $n$-th degree polynomial in the argument $z(s)$ with the
leading coefficient $\lambda_n$. Hence, it is possible to choose a
polynomial $P_n(x)$ with the property \re{2form}.

Consider now condition \re{eq_Q} for $n \to n+1$: \be A(s) \Delta
z^{n+1}(s) - C(s) \nabla z^{n+1}(s) = \lambda_{n+1} Q_{n+1}(z(s)),
\quad n=0,1,\dots \lab{eq_Q1} \ee Multiplying \re{eq_Q} by $z(s)$
and subtracting \re{eq_Q1} we get another set of conditions \be
A_1(s) z^n(s+1) + C_1(s) z^n(s-1) = R_{n+1}(z(s)),  \quad
n=0,1,\dots \lab{eq_R} \ee where $A_1(s) = A(s) \Delta z(s), \quad
C_1(s) = C(s) \nabla z(s)$, The polynomials $R_n(z)$ are $n$-th
degree polynomials $R_n(z) =\omega_n z^n + O(z^{n-1}) $, where
$\omega_n =\lambda_{n} - \lambda_{n-1}$. Note that due to
non-degeneracy condition \re{ndegn} we have $\omega_n \ne 0$ and
hence every polynomial $R_n(z)$ has exact degree $n$.

Consider first two conditions \re{eq_R} corresponding to $n=0$ and
$n=1$. These two conditions can be considered as equations for two
unknowns $A_1(s), C_1(s)$. Solving these equations we have \ba
&&A_1(s)= \frac{R_2(z(s)) -z(s-1) R_1(z(s) }{z(s+1) - z(s-1)}, \nonumber \\
&&C_1(s)=- \frac{R_2(z(s)) -z(s+1) R_1(z(s) }{z(s+1) - z(s-1)}
\lab{AC} \ea Note that these expressions are well defined for all
possible $s$ because, by non-degeneracy condition,  $z(s+1) \ne
z(s-1)$.

Hence conditions \re{eq_R} can be rewritten as \be R_2(z(s)) Y_n -
R_1(z(s)) z(s-1)z(s+1) Y_{n-1} = R_{n+1}(z(s)), \; n=2,3,\dots,
\lab{eq_R1} \ee where \be Y_n = \frac{z^n(s+1) -
z^n(s-1)}{z(s+1)-z(s-1)}. \lab{Y_n} \ee Introduce the variables
$$
u=z(s-1)z(s+1), \quad v=z(s-1) + z(s+1)
$$
Clearly $Y_n$ is a symmetric polynomial with respect to $z(s-1),
z(s+1)$ and hence it can be expressed in terms of variables $u,v$
only. Indeed, it is easily verified that $Y_n$ satisfy the
recurrence relation \be Y_{n+1} = vY_n - uY_{n-1}, \quad Y_0=0, \:
Y_1=1. \lab{rec_Y} \ee This allows us to find an explicit
expression for every $Y_n$ in terms of $u,v$. For example, $Y_2=
v, \: Y_3= v^2-u,\: Y_4= v^3 -2uv$ etc.

Return to condition \re{eq_Q}. We have explicit expressions for
coefficients $A(s),C(s)$: \ba &&A(s)=\frac{R_1 z(s-1) -
R_2}{(z(s+1)-z(s-1))(z(s+1)-z(s))}, \nonumber \\ &&C(s)=\frac{R_1
z(s+1) - R_2}{(z(s+1)-z(s-1))(z(s)-z(s-1))} \lab{AC_exp} \ea Hence
we have
$$
A(s) \Delta z(s) - C(s) \nabla z(s)= \sum_{k=0}^{n-1}
z^{n-k-1}(z(s+1)z(s-1)R_1 Y_{k-1} - R_2 Y_k) = -\sum_{k=0}^{n-1}
R_{k+1} z(s)^{n-k-1}
$$
(in the last equality we have used \re{eq_R1}). It is seen that
this expression is indeed a polynomial of degree $n$ with non-zero
leading coefficient $\lambda_n$. Thus conditions \re{eq_Q} and
\re{eq_R} are equivalent and we can use only more simple condition
\re{eq_R} for further analysis.

\section{Functional equation for the grid $z(s)$}
\setcounter{equation}{0} From \re{eq_R1} and \re{rec_Y} we find
the conditions \be R_{n+2}(z(s)) = vR_{n+1}(z(s)) - uR_n(z(s)),
\quad n=2,3,\dots \lab{rec_R} \ee These conditions form a system
of linear equations for two unknowns $u,v$. Consider the first two
equations corresponding to $n=2$ and $n=3$. There are two
possibilities:

(i) these equations are not independent. Then we should have
$R_{i+1}(x)=\tau(x) R_i(x), \; i=1,2,3$ where $\tau(x)$ is a
linear function. By induction, we then have $R_n(x) = R_1(x)
\tau^{n-1}(x),\; n=1,2,\dots$ for all $n$, where both $\tau(x)$
and $R_1(x)$ are linear functions in $x$. Now from \re{rec_R} we
have the condition \be \tau^2(z(s)) - v \tau(z(s)) +u =0
\lab{dep_con} \ee or, equivalently, \be (\tau(z(s))
-z(s+1))(\tau(z(s))-z(s-1))=0. \lab{dep_1} \ee From \re{dep_1} and
\re{AC_exp} we see that in this case either $A(s)=0$ or $C(s)=0$
for every admissible $s$. But this contradicts our non-degeneracy
assumption $A(s-1)C(s) \ne 0$. Thus the case (i) should be
excluded from consideration.

(ii) these equations are independent.  Putting $n=2,3$ in
\re{rec_R} we obtain a linear system of equations for unknowns
$u,v$ from which we find \be u=\frac{\pi_8(z(s))}{\pi_6(z(s))},
\quad v=\frac{\pi_7(z(s))}{\pi_6(z(s))}, \lab{uv_rat} \ee where
$\pi_i(x)$ are polynomials of degrees $\le i$:
$$
\pi_6= R_3^2-R_2R_4, \; \pi_7=R_4R_3 -R_2R_5, \; \pi_8 =R_3R_5 -
R_4^2
$$

Thus $u,v$ are some rational functions in the variable $z(s)$.
 In what follows we will sometimes
replace the grid $z(s)$ with independent variable $x$ (this is
possible because the grid $z(s)$ takes infinitely many different
values).

We first prove an important statement concerning possible
solutions of the system of non-linear difference equations of the
form \be z(s-1)+z(s+1)=T_1(z(s)), \quad z(s-1) z(s+1) = T_2(z(s))
\lab{sys_TT} \ee where $T_{1,2}(x)$ are some rational functions.

\begin{lem}
Assume that the system \re{sys_TT} has a solution $z(s),\;
s=s_0,s_0+1,\dots$ with infinitely many non-coinciding values
$z(s_1) \ne z(s_2)$ if $s_1 \ne s_2$. Then there are two
possibilities:

(i) either \be T_1(x) = - \frac{\alpha_2 x^2 + \alpha_4 x +
\alpha_5}{\alpha_1 x^2 + \alpha_2 x + \alpha_3}, \quad T_2(x) =
\frac{\alpha_3 x^2 + \alpha_5 x + \alpha_6}{\alpha_1 x^2 +
\alpha_2 x + \alpha_3} \lab{1_TT} \ee

with some constants $\alpha_i,i=1,\dots,6$.

In this case variables $z(s),z(s+1)$ satisfy equation \be
\Phi(z(s),z(s+1)) =0, \lab{eq_Phi} \ee where $\Phi(x,y)$ is a
non-reducible symmetric bi-quadratic polynomial: \be \Phi(x,y) =
\alpha_1 x^2y^2 + \alpha_2 xy(x+y) + \alpha_3 (x^2+y^2) + \alpha_4
xy + \alpha_5 (x+y) + \alpha_6 \lab{Phi} \ee

or

(ii) \be T_1(x) =-\frac{\alpha_2 x + \alpha_4}{\alpha_1 x +
\alpha_3} - \frac{\alpha_3 x + \alpha_4}{\alpha_1 x + \alpha_2},
\quad T_2(x) =\frac{\alpha_2 x + \alpha_4}{\alpha_1 x + \alpha_3}
\: \frac{\alpha_3 x + \alpha_4}{\alpha_1 x + \alpha_2}
\lab{TT_deg} \ee with some constants $\alpha_i, i=1,\dots,4$ such
that $\alpha_2 \ne \alpha_3$. In this case variables $z(s),z(s+1)$
satisfy equation \be \alpha_1 z(s) z(s+1) + \alpha_2 z(s) +
\alpha_3 z(s+1) + \alpha_4=0 \lab{eq_deg} \ee

\end{lem}

{\it Remark}. The case (ii) formally corresponds to a special case
of (i) when polynomial $\Phi(x,y)$ can be decomposed as a product
of two polynomials of the first degree in both variables $x,y$.

{\it Proof}. Obviously, system \re{sys_TT} is equivalent to the
statement that both $z(s+1)$ and $z(s-1)$ are roots of the
quadratic equation \be A_2(z(s)) z_{s\pm 1}^2 + A_1(z(s)) z_{s\pm
1} + A_0(z(s)) =0, \lab{quad_zzz} \ee where $A_i(x)$ are non-zero
polynomials having no common factors.

Introduce two polynomials in two variables:
$$
W_1(x,y) =A_2(x) y^2 + A_1(x) y + A_0(x), \quad W_2(x,y) \equiv
W_1(y,x) = A_2(y) x^2 + A_1(x) x + A_0(y)
$$
Equations $W_1(x,y)=0$ and $W_2(x,y)=0$ define two algebraic
curves in complex variables $x,y$. From \re{quad_zzz} it is clear
that both curves contain infinitely many common distinct points
$(x_n,y_n), n=1,2,\dots$. By the Bezout theorem this is possible
only if these curves either coincide or have a common component.

The polynomial $W_1(x,y)$ has degree 2 in variable $y$ and hence
there are two possibilities:

(i) $W_1(x,y)$ is irreducible, i.e. it cannot be decomposed into
irreducible polynomials of a lesser degree in $y$.

(ii) $W_1(x,y)$ can be presented as a product of two polynomials,
each of degree 1 in variable $y$: $W_1(x,y)=(e_1(x) y +
e_2(x))(e_3(x)y + e_4(x))$ with some polynomials $Re_i(x),
i=1,\dots,4$.

We consider these two possibilities separately. In the case (i) we
have that the polynomials $W_1(x,y)$ and $W_2(x,y)$ are both
irreducible. Hence, by the Bezout theorem, they should coincide:
\be W_1(x,y) = W_2(x,y) = W_1(y,x) \lab{sym_WW} \ee But condition
\re{sym_WW} means that the polynomial $W_1(x,y)$ is {\it
symmetric} in variables $x,y$. This is possible only if all
polynomials $A_i(x), i=0,1,2$ have degree $\le 2$ in variable $x$.
Hence, the most general expression for $W_1(x,y)$ in this case is
{\it symmetric bi-quadratic polynomial} in $x,y$: \be W_1(x,y) =
\alpha_1 x^2y^2 + \alpha_2 xy(x+y) + \alpha_3 (x^2+y^2) + \alpha_4
xy + \alpha_5 (x+y) + \alpha_6 \lab{bi_qua_W} \ee with some
constants $\alpha_i, i=1,\dots, 6$. We thus have
$$
A_2(x) = \alpha_1 x^2 + \alpha_2 x + \alpha_3, \; A_1(x) =
\alpha_2 x^2 + \alpha_4 x + \alpha_5, \; A_0(x) = \alpha_3 x^2 +
\alpha_5 x + \alpha_6
$$
and
$$
T_1(x) = -A_1(x)/A_2(x)= - \frac{\alpha_2 x^2 + \alpha_4 x +
\alpha_5}{\alpha_1 x^2 + \alpha_2 x + \alpha_3}, \quad T_2(x) =
A_0(x)/A_2(x)=\frac{\alpha_3 x^2 + \alpha_5 x + \alpha_6}{\alpha_1
x^2 + \alpha_2 x + \alpha_3}
$$
giving expression \re{1_TT}.

Consider now the case when polynomials $W_1(x,y),W_2(x,y)$ have a
nontrivial common component which doesn't coincide with both these
polynomials. Clearly, this is possible only if $W_1(x,y)$ can be
decomposed into two polynomials linear in $y$: \be W_1(x,y) =
(a_1(x)y+b_1(x))(a_2(x) y + b_2(x)) \lab{W1_fact} \ee with some
polynomials $a_{1,2}(x), b_{1,2}(x)$. By definition
$W_2(x,y)=W_1(y,x)$ and hence we have also \be W_2(x,y) =
(a_1(y)x+b_1(y))(a_2(y) x + b_2(y)) \lab{W2_fact} \ee Without loss
of generality we can assume that $a_1(x) y + b_1(x)$ is a common
component of two curves $W_1(x,y)=0$ and $W_2(x,y)=0$. Comparing
\re{W1_fact} and \re{W2_fact}, we can conclude that there are two
possibilities:

(i) either $a_1(x) y +b_1(x) = a_1(y) x + b_1(y)$;

(ii) or $a_1(x) y +b_1(x) = a_2(y) x + b_2(y)$

In case (i) we have that variables $x,y$ satisfy symmetric
polynomial relation \be \alpha_1 xy + \alpha_2(x+y) + \alpha_4 =0
\lab{sym_deg} \ee with some constants
$\alpha_1,\alpha_2,\alpha_4$. Substituting $x=z(s),y=z(s+1)$ into
\re{sym_deg} we find from \re{sym_deg} that there are only two
non-coinciding points $z(s_0)$ and $z(s_0+1)$. For all further
points we find that $z(s_0 + 2j) = z(s_0)$ and
$z(s_0+2j+1)=z(s_0+1)$ for all integer $j$. But this contradicts
our assumption that there are infinitely many distinct points
belonging to the curves. Thus the case (i) is impossible.

In the case (ii) the polynomials $a_{1,2}(x), b_{1,2}(x)$ should
be linear in $x$ and we have that $z(s),z(s+1)$ satisfy the
relation \be \alpha_1 z(s) z(s+1)+ \alpha_2 z(s) + \alpha_3 z(s+1)
+ \alpha_4 \lab{nsym_deg} \ee where $\alpha_3 \ne \alpha_2$ in
order to prevent impossible case (i). This case corresponds to
\re{TT_deg} and \re{eq_deg}. Thus the Lemma is proven.

It is interesting to find explicit solutions in both cases (i) and
(ii) of the Lemma. The case (i) corresponds to a parametrization
of symmetric Euler-Baxter bi-quadratic curve $\Phi(z(s),z(s+1))=0$
with $\Phi(x,y)$ given by \re{Phi}. This problem was already
solved by Baxter \cite{Bax} in his famous solution of the 8-vertex
model. Explicitly \be z(s) = \kappa \: \phi(\beta_1 s + \beta_0)
\lab{Bax} \ee with some parameters $\kappa, \beta_1, \beta_0$.
Here $\phi(z)$ is an even elliptic function of the second order
(i.e. having exactly two poles in the fundamental parallelogram).
Recall that (up to an arbitrary factor) any even elliptic function
of the second order can be presented in the form \cite{WW} \be
\phi(z) = \frac{\sigma(z-e)\sigma(z+e)}{\sigma(z-d) \sigma(z+d)}
\lab{2_degree} \ee Recently it was shown that the elliptic grid
$z(s)$ described by \re{Bax} appears naturally in theory of
biorthogonal rational functions with the duality property
\cite{SZ2}, \cite{SZ_Ky}. Note that in a special case
$\alpha_1=\alpha_2=0$ the rational functions $T_1(x),T_2(x)$
become linear and quadratic polynomials. In this case solution for
$z(s)$ is expressed in terms of elementary functions (see below).

For the case (ii) of the Lemma the solution can be easily found in
terms of elementary functions of $s$ (we will  not describe these
solutions in details because they can be obtained from the
elliptic solutions by a limiting procedure).

Now return to condition \re{rec_R} and consider first the case (i)
of the Lemma. We can rewrite \re{rec_R} in the form \be 1=v(x)
R_{n+1}(x)/R_{n+2}(x)-u(x) R_n(x)/R_{n+2}(x), \lab{rec_R1} \ee
where
$$
v(x) =-\frac{\alpha_2 x^2 + \alpha_4 x + \alpha_5}{\alpha_1 x^2 +
\alpha_2 x + \alpha_3}, \quad u(x) = \frac{\alpha_3 x^2 + \alpha_5
x + \alpha_6}{\alpha_1 x^2 + \alpha_2 x + \alpha_3}
$$
Assume first that $\alpha_1 \ne 0$. Then for $x \to \infty$ it is
seen that rhs of  \re{rec_R1} tends to 0 which contradicts to lhs
of \re{rec_R1}. Thus necessarily $\alpha_1 =0$. Assume now that
$\alpha_1 = 0$ and $\alpha_2 \ne 0$. Then again for $x \to \infty$
we obtain from \re{rec_R1} the condition (recall that $R_n(x) =
\omega_n x^n + O(x^{n-1})$, where $\omega_n = \lambda_n -
\lambda_{n-1}$)
$$
1=-\omega_{n+1}/\omega_{n+2}
$$
whence $\lambda_n = \lambda_{n+2}$ for all $n=2,3,\dots$. But this
contradicts our condition of non-degeneracy of the spectrum
$\lambda_n$. We thus have necessarily $\alpha_1=\alpha_2 =0$. But
in this case $v(x), u(x)$ become polynomials of the first and
second degrees: \be v(x) = -\xi x - \eta, \quad u(x)= x^2 + \eta x
+ \zeta \lab{vu_pol} \ee where $\xi=\alpha_4 /\alpha_3, \: \eta =
\alpha_5 /\alpha_3, \: \zeta = \alpha_6 /\alpha_3 $ and equations
for the grid become \be z(s-1) + z(s+1) = -\xi z(s) - \eta, \quad
z(s-1) z(s+1) = z^2(s) + \eta z(s) + \zeta \lab{NSU_rec} \ee with
arbitrary complex parameters $\xi,\eta,\zeta$. Equivalently,
variables $z(s), z(s+1)$ belong to a non-degenerating conic (i.e.
ellipsis, hyperbola or parabola): \be z^2(s+1) + z^2(s)
+\eta(z(s+1) + z(s)) + \xi z(s) z(s+1) + \zeta =0 \lab{Mag_rec}
\ee which is symmetric with respect to $z(s), z(s+1)$ (this means
that the plot of this conic in Cartesian co-ordinates
$x=z(s),y=z(s+1)$ is symmetric with respect to the line $y=x$).
Equations \re{NSU_rec} and \re{Mag_rec} were studied in
\cite{NSU}, \cite{Mag1}, \cite{Mag3}. In these works it was shown
that all nondegenerate solutions of these equations can be
presented in  the form \be z(s) = C_1 q^s + C_2 q^{-s} + C_0
\lab{z_trig} \ee or \be z(s) = C_2 s^2 + C_1 s + C_0 \lab{z_quad}
\ee  or \be z(s) = (-1)^s \: (C_2 s^2 + C_1 s + C_0) \lab{z_quad1}
\ee   with some constants $C_0,C_1,C_2$. The first case
\re{z_trig} occurs if $\xi = q+q^{-1}$, where $q\ne \pm 1$ (i.e.
$\xi \ne \pm 2$). The second case \re{z_quad} occurs if $\xi=-2$
and the third case \re{z_quad1} occurs if $\xi =2$. All these
cases exhaust possible types of the  Askey-Wilson grids
\cite{Mag3}.

Note that when $C_1C_2 =0$ in \re{z_trig} we obtain so-called
exponential grids, say $z(s) = C_1 q^s + C_0$. Similarly, when
$C_2=0$ in \re{z_quad} or \re{z_quad1} we obtain the linear grid:
$z(s) = C_1 s + C_0$ or $z(s)=(-1)^s (C_1 s + C_0)$. However, in
these case the conic \re{Mag_rec} becomes degenerated - it divided
into two lines. This corresponds to the case (ii) of the Lemma
(see below).

Now substituting $v(x),u(x)$ into \re{rec_R} we obtain that for
arbitrary given polynomials $R_1(x),R_2(x)$ one can construct
uniquely the polynomial $R_n(x) = \omega_n x^n + O(x^{n-1})$ with
leading coefficient satisfying the recurrence relation (which
easily follows from \re{rec_R} for given $v(x)u(x)$): \be
\omega_{n+2} + \xi \omega_{n+1} + \omega_n =0, \quad n=2,3,\dots
\lab{omega_rec} \ee General solution for $\omega_n$ can be easily
found from \re{omega_rec}: if $\xi = q+ q^{1} \ne \pm 2$ we have
\be \omega_n = G_1 q^n + G_2 q^{-n} \lab{om_trig} \ee with
arbitrary $G_1,G_2$. If $\xi =-2$ then \be \omega_n = G_1 n + G_0
\lab{om_lin} \ee and if $\xi=2$ then \be \omega_n = (-1)^n(G_1 n +
G_0) \lab{om_lin1} \ee From $\omega_n=\lambda_n-\lambda_{n-1}$ we
can easily reconstruct the spectrum $\lambda_n$ which has the same
functional dependence on $n$ as the Askey-Wilson grid $z(s)$ has
on $s$.

Finally, we should consider the case (ii) of the Lemma. In this
case similar considerations lead to conclusion that $\alpha_1=0$.
Then condition \re{eq_deg} becomes \be \alpha_2 z(s) + \alpha_3
z(s+1) + \alpha_4=0 \lab{exp_eq} \ee which describes exponential
or linear grids $z(s)$. Thus the case (ii) can be considered as a
degeneration of the case (i). It should be noted that linear and
exponential grids are described by {\it non-symmetric} recurrence
relations \re{exp_eq} (with respect to $z(s),z(s+1)$). The reason
is that in this case the conic \re{Mag_rec} is degenerated to a
two straight lines each of which is non-symmetric .

We see that in both cases (i) and (ii) of the Lemma solutions
$R_n(x)$ of the recurrence relation \re{rec_R} are indeed
polynomials of exact degree $n$. Hence, by previous
considerations, we obtain a unique set of polynomials $P_n(x)$
which are solutions of equation \re{2form}. From explicit form
\re{AC_exp} of the coefficients $A(s),C(s)$ we can conclude that
they coincide with those defining the Askey-Wilson polynomials
\cite{NSU}, \cite{Mag1}, \cite{Mag3}.

We thus proved that under some non-degeneracy conditions, the only
admissible grid is the Askey-Wilson grid and corresponding
polynomials $P_n(x)$ coincide with the Askey-Wilson polynomials.

\section{Concluding remarks}
\setcounter{equation}{0} The authors of \cite{NSU} exploited
relations \re{2form} as a starting point in their approach to
construction of the Askey-Wilson polynomials. In a slightly
different manner, Magnus in \cite{Mag1}, \cite{Mag3} derived
relations \re{Mag_rec} from the following requirement: find all
the grids $z(s)$ and $y(s)$, such that for {\it any} polynomial
$P_n(x)$ of degree $n$ we have the property \be \frac{P_n(z(s+1))-
P_n(z(s))}{z(s+1)-z(s)} = T_{n-1}(y(s)), \lab{Mag_PQ} \ee where
$T_{n-1}(x)$ is a polynomial of degree $n-1$ and relation
\re{Mag_PQ} should be valid for all $n=1,2,\dots$ and for
infinitely many distinct values $s$ of the grids $z(s)$ and
$y(s)$. Relation \re{Mag_PQ} can be also presented in the form \be
{\cal D}_s P_n(x) = T_{n-1}(y(s)) \lab{DPQ} \ee where ${\cal D}_s$
stands for "discrete derivation" operator which acts on the space
of function $f(x)$ as
$$
{\cal D}_s f(x) \equiv  \frac{f(z(s+1))- f(z(s))}{z(s+1) - z(s)}
$$
For the AW-grid it was known that the operator ${\cal D}_s$
satisfies property \re{DPQ}. Magnus proved that these grids are
{\it the only} preserving property \re{DPQ}.

On the other hand, it was noted in \cite{NSU} that if polynomials
$P_n(x)$ satisfy the AW-equation \re{2form} then the new
polynomials $T_n(x)$ obtained from $P_n(x)$ by \re{Mag_PQ} also
satisfy AW-equation \re{2form} but with different coefficients
$A(s),C(s)$. This property can be considered as a covariance of
the Askey-Wilson equation \re{2form} with respect to the discrete
Darboux transformation (see, e.g. \cite{SZ_MAA}).

However the property \re{Mag_PQ} cannot be directly derived from
equation \re{2form} if the grid $z(s)$ is not concretized. This is
why derivation of the necessity of the AW-grid for equation
\re{2form} is not quite elementary and needs rather involved
technique which was demonstrated in the present paper.

\bigskip
\section{Acknowledgments}

The authors thank the referees for their critical readings of this paper, which led to its improvement. A.Zh. thanks Centre de recherches math\'ematiques of the
Universit\'e de Montr\'eal for hospitality.


\bb{99}

\bi{Al} W.A. Al-Salam, {\it Characterization theorems for
orthogonal polynomials}, in: P.Nevai (ed.), "Orthogonal
Polynomials: Theory and Practice", NATO ASI Series C: Mathematical
and Physical Sciences, vol. 294. Kluwer Academic Publishers,
pp.1-24.

\bi{AW} R.~Askey and J.~Wilson, {\it Some basic hypergeometric orthogonal
polynomials that generalize Jacobi polynomials}, Mem. Amer. Math. Soc. {\bf
54}, No. 319, (1985), 1-55.

\bi{BI} E. Bannai and T. Ito, {\it Algebraic Combinatorics I:
Association Schemes} (Benjamin \& Cummings, 1984).

\bi{Bax} R.~Baxter, {\it Exactly Solvable Models in Statistical
Mechanics}. Academic Press, 1982.

\bi{Bo} S.~Bochner, {\it \"Uber Sturm-Liouvillesche Polynomsysteme}, Math.
Zeit., {\bf 29} (1929), 730-736.

\bi{Chi} T.Chihara, An Introduction
to Orthogonal Polynomials, (Gordon and
Breach, 1978).



\bi{GH} F.Alberto Gr\"unbaum and Luc Haine, {\it The q-version of
a theorem of Bochner}, J. Comput. Appl. Math. {\bf 68} (1996),
103-114.


\bi{Ismail} M.E.H.Ismail, {\it A generalization of a theorem of
Bochner}, J. Comput. Appl. Math. {\bf 159} (2003) 319–324.

\bi{KS} Koekoek R and Swarttouw R F 1994 {\it The Askey scheme of
hypergeometric orthogonal \\polynomials and its q-analogue}, Report 94-05,
Faculty of Technical Mathematics and Informatics, Delft University of
technology.

\bi{Leo} D.A.Leonard, {\it Orthogonal polynomials, duality and association
schemes}, SIAM J.Math.Anal. {\bf 13} (1982), 656-663.

\bi{Mag1} A.Magnus, {\it Associated Askey-Wilson polynomials as Laguerre-Hahn
orthogonal polynomials}, pp 261-278, in: Orthogonal Polynomials and their
applications (ed. by M.Alfaro et al., 1988); Lect.Notes in Mathem. {\bf 1329},
(Springer, Berlin, 1988).

\bi{Mag3} A.P. Magnus, {\it Painlev\'e-type differential equations
for the recurrence coefficients of semi-classical orthogonal
polynomials}, J. Comp. Appl. Math. {\bf 57} (1995), 215-237.


\bi{NSU} A.F.Nikiforov, S.K.Suslov, and V.B.Uvarov, {\it Classical Orthogonal
Polynomials of a Discrete Variable}, Springer, 1991.

\bi{SZ_MAA} V.Spiridonov and A. Zhedanov, {\it Discrete Darboux
transformations, the discrete-time Toda lattice, and the
Askey-Wilson polynomials}. Methods Appl. Anal. {\bf 2} (1995),
369--398.

\bibitem{SZ2} V.P. Spiridonov and A.S. Zhedanov, {\it Generalized eigenvalue problem
and a new family of rational functions biorthogonal on elliptic
grids}, Special Functions 2000, Kluwer, Dordrecht (2001), pp.
365--388.

\bibitem{SZ_Ky} V.P. Spiridonov and A.S. Zhedanov, {\it To the theory of biorthogonal
rational functions}, RIMS Kokyuroku {\bf 1302} (2003), 172--192.


\bi{Terw} P.Terwilliger, {\it Two linear transformations each
tridiagonal with respect to an eigenbasis of the other}. Linear
Algebra Appl. 330 (2001), no. 1-3, 149--203.

\bi{Terw1} P.Terwilliger, {\it Two linear transformations each
tridiagonal with respect to an eigenbasis of the other: comments
on the split decomposition}, J.Comput.Appl.Math., {\bf 178}
(2005), 437--452.

\bi{WW} E.~T.~Whittacker and G.~N.~Watson, {\it A Course of Modern
Analysis}, Fourth Edition, Cambridge, University Press, 1927.

\bi{Wil} J.A. Wilson, {\it Some hypergeometric orthogonal
polynomials}. SIAM J. Math.Anal., {\bf 11} (1980), 690-701.

\eb

\end{document}